\documentclass[12pt]{article}
\usepackage{amsfonts,amscd,a4}
\usepackage{color}
\usepackage{graphicx}
\usepackage{theorem}
\newtheorem{theo}{Theorem}
\newtheorem{prop}[theo]{Proposition}
\newtheorem{lemma}[theo]{Lemma}
\newtheorem{coro}[theo]{Corollary}

{\theorembodyfont{\rm}

}

\newcommand{\cA}{{\mathcal A}}
\newcommand{\cB}{{\mathcal B}}
\newcommand{\cC}{{\mathcal C}}

\newcommand{\cF}{{\mathcal F}}
\newcommand{\cG}{{\mathcal G}}

\newcommand{\cJ}{{\mathcal J}}

\newcommand{\cS}{{\mathcal S}}

\newcommand{\eT}{{\sf T}}

\newcommand{\sC}{{\mathbb C}}

\newcommand{\sN}{{\mathbb N}}

\newcommand{\sT}{{\mathbb T}}

\newcommand{\sZ}{{\mathbb Z}}
\newcommand{\qed}{\rule{1ex}{1ex}}
\newcommand{\alg}{\mbox{\rm alg} \,}

\newcommand{\im}{\mbox{\rm im} \,}

\begin{document}
\title{The universal algebra generated by a power partial isometry}
\author{Steffen Roch}
\date{Dedicated to Albrecht B\"ottcher on the occassion of his 60th birthday.}
\maketitle
\begin{abstract}
A power partial isometry (PPI) is an element $v$ of a $C^*$-algebra with the property that every power $v^n$ is a partial isometry. The goal of this paper is to identify the universal $C^*$-algebra generated by a PPI with (a slight modification of) the algebra of the finite sections method for Toeplitz operators with continuous generating function, as first described by Albrecht B\"ottcher and Bernd Silbermann in \cite{BSi0}.  
\end{abstract}
{\bf 2010 AMS-MSC:} 46L05, 47B35, 65R20
\section{Introduction}
Let $\cA$ be a $C^*$-algebra. An element $V$ of $\cA$ is called a partial isometry if $vv^*v = v$. Simple examples show that a power of a partial isometry needs not to be a partial isometry again. One therefore calls $v$ a {\em power partial isometry} (PPI) if every power of $v$ is a partial isometry again. \\[3mm]
{\bf Examples.} $(a)$ In a $C^*$-algebra with identity element $e$, every unitary element $u$ (i.e. $u^*u = uu^* = e$) is a PPI. In particular, the function $u: t \mapsto t$ is a unitary element of the algebra $C(\sT)$ of the continuous functions on the complex unit circle $\sT$, and the operator $U$ of multiplication by the function $u$ is a unitary operator on the Hilbert space $L^2(\sT)$ of the squared integrable functions on $\sT$. \\[1mm]
$(b)$ In a $C^*$-algebra with identity element $e$, every isometry $v$ (i.e. $v^*v = e$) and every co-isometry $v$ (i.e. $vv^* = e$) is a PPI. In particular, the operators $V : (x_0, \, x_1, \, \ldots) \mapsto (0, \, x_0, \, x_1, \, \ldots)$ and $V^* : (x_0, \, x_1, \, x_2, \, \ldots) \mapsto (x_1, \, x_2, \, \ldots)$ of forward and backward shift, respectively, are PPIs on the Hilbert space $l^2(\sZ^+)$ of the squared summable sequences on the non-negative integers. \\[1mm]
$(c)$ The matrix $V_n := (a_{ij})$ with $a_{i+1,i} = 1$ and $a_{ij} = 0$ if $i \neq j+1$, considered as an element of the algebra $\sC^{n \times n}$ of the complex $n \times n$ matrices, is a PPI. \\[1mm]
$(d)$ If $v_i$ is a PPI in a $C^*$-algebra $\cA_i$ for every $i$ in an index set $I$, then $(v_i)_{i \in I}$ is a PPI in the direct product $\prod_{i \in I} \cA_i$. In particular, the operator $(V, \, V^*)$, considered as an element of $L(l^2(\sZ^+)) \times L(l^2(\sZ^+))$, is a PPI. \\[1mm]
Note that the PPI $V_n$ in $(c)$ and $(V, \, V^*)$ in $(d)$ are neither isometric nor co-isometric.  \hfill \qed \\[3mm]
The goal of the present paper is to describe the universal $C^*$-algebra generated by a PPI. Recall that a $C^*$-algebra $\cA$ generated by a PPI $v$ is {\em universal} if, for every other $C^*$-algebra $\cB$ generated by a PPI $w$, there is a $^*$-homomorphism from $\cA$ to $\cB$ which sends $v$ to $w$. The universal algebra generated by a unitary resp. isometric element is defined in an analogous way. The existence of a universal algebra generated by a PPI is basically a consequence of Example $(d)$.

It follows from the Gelfand-Naimark theorem that the universal algebra generated by a unitary element is $^*$-isomorphic to the algebra $C(\sT)$, generated by the unitary function $u$. Coburn \cite{Cob1} identified the universal algebra generated by an isometry as the Toeplitz algebra $\eT(C)$ which is the smallest $C^*$-subalgebra of $L(l^2(\sZ^+))$ which contains the isometry $V$, the shift operator. This algebra bears its name since it can be described as the smallest $C^*$-subalgebra of $L(l^2(\sZ^+))$ which contains all Toeplitz operators $T(a)$ with generating function $a \in C(\sT)$. Recall that the Toeplitz operator with generating function $a \in L^1(\sT)$ is given the matrix $(a_{i-j})_{i,j = 0}^\infty$ where $a_k$ stands for the $k$th Fourier coefficient of $a$. This operator is bounded on $l^2(\sZ^+)$ if and only if $a \in L^\infty(\sT)$ (see \cite{BSi1,BSi2}).

We will see that the universal algebra of a PPI is also related with Toeplitz operators, via the finite sections discretization with respect to the sequence of the projections $P_n : (x_0, \, x_1, \, \ldots) \mapsto (x_0, \, \ldots, \, x_{n-1}, 0, \, 0, \, \ldots)$ on $l^2(\sZ^+)$. Write $\cF$ for the set of all bounded sequences $(A_n)_{n \ge 1}$ of operators $A_n \in L(\im P_n)$ and $\cG$ for the set of all sequences $(A_n) \in \cF$ with $\|A_n\| \to 0$. Provided with entry-wise defined operations and the supremum norm, $\cF$ becomes a $C^*$-algebra and $\cG$ a closed ideal of $\cF$. Since $L(\im P_n)$ is isomorphic to $\sC^{n \times n}$, we can identify $\cF$ with the direct product and $\cG$ with the direct sum of the algebras $\sC^{n \times n}$ for $n \ge 1$. Now consider the smallest $C^*$-subalgebra $\cS(\eT(C))$ of $\cF$ which contains all sequences $(P_n T(a) P_n)$ with $a \in C(\sT)$ and its $C^*$-subalgebra $\cS_{\ge 2}(\eT(C))$ which is generated by the sequence $(P_n V P_n)$ (note that $V$ is the Toeplitz operator with generating function $t \mapsto t$). With these notations, the main result of the present paper can be formulated as follows.
\begin{theo} \label{t1}
The universal algebra generated by a PPI is $^*$-isomorphic to the $C^*$-algebra $\cS_{\ge 2}(\eT(C))$ generated by the PPI $(P_n V P_n)$.
\end{theo}
For a general account on $C^*$-algebras generated by partial isometries, with special emphasis on their relation to graph theory, see \cite{ChJ1}.
 
Before going into the details of the proof of Theorem \ref{t1}, we provide some basic (and well known) facts on the algebras $\cS (\eT(C))$ and $\cS_{\ge 2}(\eT(C))$. Since the first entry of the sequence $(P_n V P_n)$ is zero, the first entry of every sequence in $\cS_{\ge 2}(\eT(C))$ is zero. So we can omit the first entry and consider the elements of $\cS_{\ge 2}(\eT(C))$ as sequences labeled by $n \ge 2$ (whence the notation). In fact this is the only difference between the algebras $\cS (\eT(C))$ and $\cS_{\ge 2}(\eT(C))$.
\begin{prop} \label{p2}
$\cS_{\ge 2}(\eT(C))$ consists of all sequences $(A_n)_{n \ge 2}$ where $(A_n)_{n \ge 1}$ is a sequence in $\cS(\eT(C))$.
\end{prop}
The sequences in $\cS(\eT(C))$ are completely described in the following theorem, where we let $R_n$ denote the operator $(x_0, \, x_1, \, \ldots) \mapsto (x_{n-1}, \, \ldots, \, x_0, 0, \, 0, \, \ldots)$ on $l^2(\sZ^+)$. Further we set $\tilde{a}(t) := a(t^{-1})$ for every function $a$ on $\sT$. This description was found by A. B\"ottcher and B. Silbermann and first published in their 1983 paper \cite{BSi0} on the convergence of the finite sections method for quarter plane Toeplitz operators (see also \cite{HRS2}, Section 1.4.2).
\begin{prop} \label{p3}
The algebra $\cS(\eT(C))$ consists of all sequences $(A_n)_{n \ge 1}$ of the form
\begin{equation} \label{e1}
(A_n) = (P_n T(a) P_n + P_n K P_n + R_n L R_n + G_n)
\end{equation}
where $a \in C(\sT)$, $K$ and $L$ are compact operators, and $(G_n) \in \cG$. The representation of a sequence $(A_n) \in \cS(\eT(C))$ in this form is unique.
\end{prop}
\begin{coro} \label{c4}
$\cG$ is a closed ideal of $\cS(\eT(C))$, and the quotient algebra $\cS(\eT(C))/\cG$ is $^*$-isomorphic to the $C^*$-algebra of all pairs
\begin{equation} \label{e2}
(T(a) + K, \, T(\tilde{a}) + L) \in L(l^2(\sZ^+)) \times L(l^2(\sZ^+))
\end{equation}
with $a \in C(\sT)$ and $K, \, L$ compact. In particular, the mapping which sends the sequence $(\ref{e1})$ to the pair $(\ref{e2})$ is a $^*$-homomorphism from $\cS(\eT(C))$ onto $\cS(\eT(C))/\cG$ with kernel $\cG$.
\end{coro}
It is not hard to see that the algebra of all pairs (\ref{e2}) is just the smallest $C^*$-subalgebra of $L(l^2(\sZ^+)) \times L(l^2(\sZ^+))$ that contains the PPI $(V, \, V^*)$.
\begin{coro} \label{c5}
The set $\cJ$ of all pairs $(K, \, L)$ of compact operators $K, \, L$ is a closed ideal of $\cS(\eT(C))/\cG$. The quotient algebra $(\cS(\eT(C))/\cG)/\cJ$ is $^*$-isomorphic to $C(\sT)$. In particular, the mapping which sends the pair $(\ref{e2})$ to the function $a$ is a $^*$-homomorphism from $\cS(\eT(C))/\cG$ onto $C(\sT)$ with kernel $\cJ$.
\end{coro}
Observe that all of the above examples $(a)$ - $(d)$ appear somewhere in the algebra $\cS(\eT(C))$ and its quotients.
\section{Elementary properties of PPI}
Our first goal is a condition ensuring that the product of two partial isometries is a partial isometry again.
\begin{prop} \label{p6}
Let $u, \, v$ be partial isometries. Then $uv$ is a partial isometry if and only if
\begin{equation} \label{e3}
u^*u vv^* = vv^* u^*u,
\end{equation}
i.e. if the initial projection $u^*u$ of $u$ and the range projection $vv^*$ of $v$ commute.
\end{prop}
{\bf Proof.} Condition (\ref{e3}) implies that
\[
(uv)(uv)^*(uv) = uvv^*u^*uv = (uu^*u) (vv^*v) = uv;
\]
hence, $uv$ is a partial isometry. Conversely, if $uv$ is a partial isometry, then a simple calculation gives 
\[
v^* (vv^*u^*u - u^*uvv^*)(u^*uvv^* - vv^*u^*u) v = 0.
\]
With the $C^*$-axiom we conclude that $v^* (vv^*u^*u - u^*uvv^*) = 0$, hence $vv^* (vv^*u^*u - u^*uvv^*) = 0$, which finally gives 
\[
vv^*u^*u = vv^* u^*u vv^*.
\]
The right-hand side of this equality is selfadjoint; so must be the left-hand side. Thus, $vv^*u^*u = (vv^*u^*u)^* = u^*u vv^*$, which is condition (\ref{e3}). \hfill \qed \\[3mm]
%
%
In particular, if $v$ is a partial isometry, then $v^2$ is a partial isometry if and only if
\begin{equation} \label{e4}
v^*v vv^* = vv^* v^*v.
\end{equation}
\begin{prop} \label{p8}
Let $v$ be a partial isometry with property $(\ref{e4})$ (e.g. a PPI). Then
\[
e := v^*v + vv^* - v^*vvv^* = v^*v + vv^* - vv^* v^*v
\]
is the identity element of the $C^*$-algebra generated by $v$. Moreover, 
\[
p := vv^* - vv^*v^*v = e - v^*v \quad \mbox{and} \quad \tilde{p} := v^*v - v^*vvv^* = e - vv^*
\] 
are mutually orthogonal projections (meaning that $p \tilde{p} = \tilde{p} p = 0$).
\end{prop}
{\bf Proof.} Condition (\ref{e4}) implies that $e$ is selfadjoint. Further,
\[
ve = vv^*v + vvv^* - vv^*vvv^* = v + vvv^* - vvv^* = v
\]
and, similarly, $v^*e = v^*$. Taking adjoints it follows that $ev^* = v^*$ and $ev = v$, and $e$ is the identity element. The remaining assertions are also easy to check. \hfill \qed \\[3mm]
We will often use the notation $v^{*n}$ instead of $(v^*)^n$.
\begin{prop} \label{p9}
$(a)$ If $v$ is a PPI, then
\begin{equation} \label{e5}
v^{*k} v^k v^n v^{*n} = v^n v^{*n} v^{*k} v^k \quad \mbox{for all} \; k, \, n \ge 1.
\end{equation}
$(b)$ If $v$ is a partial isometry and if $(\ref{e5})$ holds for $k=1$ and for every $n \ge 1$, then $v$ is a PPI.
\end{prop}
{\bf Proof.} Assertion $(a)$ is a consequence of Proposition \ref{p6} (the partial isometry $v^{n+k}$ is the product of the partial isometries $v^k$ and $v^n$). Assertion $(b)$ follows easily by induction. For $k=1$, condition (\ref{e5}) reduces to
\[
(v^*v) (v^n v^{*n}) = (v^n v^{*n}) (v^*v).
\]
Thus if $v$ and $v^n$ are partial isometries, then $v^{n+1}$ is a partial isometry by Proposition \ref{p6}. \hfill \qed 
\begin{lemma} \label{l10}
If $v$ is a PPI, then $(v^n v^{*n})_{n \ge 0}$ and $(v^{*n} v^n)_{n \ge 0}$ are decreasing sequences of pairwise commuting projections.
\end{lemma}
{\bf Proof.} The PPI property implies that the $v^n v^{*n}$ are projections and that
\[
v^n v^{*n} n^{n+k} (v^*)^{n+k} = (v^n v^{*n} v^n) v^k (v^*)^{n+k} = v^n v^k (v^*)^{n+k} = n^{n+k} (v^*)^{n+k}
\]
for $k, \, n \ge 0$. The assertions for the second sequence follow similarly. \hfill \qed
\section{A distinguished ideal}
Let $\cA$ be a $C^*$-algebra generated by a PPI $v$. By $\alg(v, \, v^*)$ we denote the smallest (symmetric, not necessarily closed) subalgebra of $\cA$ which contains $v$ and $v^*$. Further we write $\sN_v$ for the set of all non-negative integers such that $p v^n \tilde{p} \neq 0$. From Proposition \ref{p8} we know that $0 \not \in \sN_v$. Finally, we set
\[
\pi_n := p v^n \tilde{p} v^{*n} p \quad \mbox{and} \quad \tilde{\pi}_n := \tilde{p} v^{*n} p v^n \tilde{p}.
\]
\begin{prop} \label{p11}
$(a)$ The element $p v^n \tilde{p}$ is a partial isometry with initial projection $\tilde{\pi}_n$ and range projection $\pi_n$. Thus, the projections $\pi_n$ and $\tilde{\pi}_n$ are Murray-von Neumann equivalent in $\cA$, and they generate the same ideal of $\cA$. \\[1mm]
$(b)$ $\pi_m \pi_n = 0$ and $\tilde{\pi}_m \tilde{\pi}_n = 0$ whenever $m \neq n$.
\end{prop}
{\bf Proof.} $(a)$ By definition,
\[
\pi_n = p v^n \tilde{p} v^{*n} p = p v^n (e - vv^*) v^{*n} p = p (v^n v^{*n} - v^{n+1} (v^*)^{n+1}) p.
\]
Since $p = e - vv^*$ and $v^n v^{*n}$ commute by Proposition \ref{p9},
\[
\pi_n = p (v^n v^{*n} - v^{n+1} (v^*)^{n+1}) = (v^n v^{*n} - v^{n+1} (v^*)^{n+1}) p.
\]
Being a product of commuting projections (Lemma \ref{l10}), $\pi_n$ is itself a projection. Analogously, $\tilde{\pi}_n$ is a projection. Thus, $p v^n \tilde{p}$ is a partial isometry, and $\pi_n$ and $\tilde{\pi}_n$ are Murray-von Neumann equivalent. Finally, the equality
\[
\pi_n = \pi_n^2 = (p v^n \tilde{p} v^{*n} p)^2 = p v^n \tilde{\pi}_n v^{*n} p
\]
shows that $\pi_n$ belongs to the ideal generated by $\tilde{\pi}_n$. The reverse inclusion follows analogously. Assertion $(b)$ is again a simple consequence of Lemma \ref{l10}.   \hfill \qed \\[3mm]
Let $\cC_n$ denote the smallest closed ideal of $\cA$ which contains the projection $\pi_n$ (likewise, the projection $\tilde{\pi}_n$). We want to show that $\cC_n$ is isomorphic to $\sC^{(n+1) \times (n+1)}$ whenever $n \in \sN_v$ (Proposition  \ref{p18} below). For we need to establish a couple of facts on (finite) words in $\alg(v, \, v^*)$.
\begin{lemma} \label{l12}
Let $a, \, b, \, c$ be non-negative integers. Then
\[
v^{*a} v^b v^{*c} = \left\{ \begin{array}{lcl}
(v^*)^{a-b+c}  & \mbox{if} & \min \{a, \, c\} \ge b, \\
v^{b-a} v^{*c} & \mbox{if} & a \le b \le c, \\
v^{*a} v^{b-c} & \mbox{if} & a \ge b \ge c
\end{array} \right.
\]
and
\[
v^a v^{*b} v^c = \left\{ \begin{array}{lcl}
v^{a-b+c}       & \mbox{if} & \min \{a, \, c\} \ge b, \\
v^a (v^*)^{b-c} & \mbox{if} & a \ge b \ge c, \\
(v^*)^{b-a} v^c & \mbox{if} & a \le b \le c.
\end{array} \right.
\]
\end{lemma}
{\bf Proof.} Let $\min \{a, \, c\} \ge b$. Then
\[
v^{*a} v^b v^{*c} = (v^*)^{a-b} v^{*b} v^b v^{*b} (v^*)^{c-b} = (v^*)^{a-b} v^{*b} (v^*)^{c-b} = (v^*)^{a-b+c},
\]
where we used that $v^{*b}$ is a partial isometry. If $a \le b \le c$, then
\[
v^{*a} v^b v^{*c} = v^{*a} v^a v^{b-a} (v^*)^{b-a} (v^*)^{c-b+a} = v^{b-a} (v^*)^{b-a} v^{*a} v^a (v^*)^{c-b+a}
\]
by Proposition \ref{p9} $(a)$. Thus,
\[
v^{*a} v^b v^{*c} = v^{b-a} (v^*)^{b-a} v^{*a} v^a v^{*a} (v^*)^{c-b} = v^{b-a} (v^*)^{b-a} v^{*a} (v^*)^{c-b} = v^{b-a} v^{*c}.
\]
Similarly, $v^{*a} v^b v^{*c} = v^{*a} v^{b-c}$ if $a \ge b \ge c$. The second assertion of the lemma follows by taking adjoints. \hfill \qed \\[3mm]
Every word in $\alg(v, \, v^*)$ is a product of powers $v^n$ and $v^{*m}$. Every product $v^a v^{*b} v^c$ and $v^{*a} v^b v^{*c}$ of three powers can be written as a product of at most two powers if one of the conditions
\begin{equation} \label{e6}
\min \{a, \, c\} \ge b \quad \mbox{or} \quad a \le b \le c \quad \mbox{or} \quad a \ge b \ge c
\end{equation}
in Lemma \ref{l12} is satisfied. Since (\ref{e6}) is equivalent to $\max \{a, \, c\} \ge b$, such a product can not be written as a product of less than three powers by means of Lemma \ref{l12} if $\max \{a, \, c\} < b$. Since it is not possible in a product $v^a v^{*b} v^c v^{*d}$ or $v^{*a} v^b v^{*c} v^d$ of four powers that $\max \{a, \, c\} < b$ {\em and} $\max \{b, \, d\} < c$, one can shorten every product of powers $v^n$ and $v^{*m}$ to a product of at most three powers. Summarizing we get the following lemma.
\begin{lemma} \label{l13}
Every finite word in $\alg (v, \, v^*)$ is of the form $v^a v^{*b}$ or $v^{*b} v^a$ with $a, \, b \ge 0$ or of the form $v^a v^{*b} v^c$ or $v^{*a} v^b v^{*c}$ with $0 < \min \{a, \, c\} \le \max \{a, \, c\} < b$.
\end{lemma}
\begin{coro} \label{c14}
Let $w$ be a word in $\alg (v, \, v^*)$. \\[1mm]
$(a)$ If $pwp \neq 0$, then $w = v^a v^{*a}$ for some $a \ge 0$. \\
$(b)$ If $\tilde{p} w \tilde{w} \neq 0$, then $w = v^{*a} v^a$ with some $a \ge 0$.
\end{coro}
{\bf Proof.} We only check assertion $(a)$. By the preceding lemma, $w$ is a product of at most three powers $v^a v^{*b} v^c$ or $v^{*a} v^b v^{*c}$. First let $w = v^a v^{*b} v^c$.
Since $vp = pv^* = 0$, we conclude that $c=0$ if $pwp \neq 0$. Writing
\[
pwp = \left\{ \begin{array}{lcl}
p v^a v^{*a} (v^*)^{b-a} p = v^a v^{*a} p (v^*)^{b-a} p & \mbox{if} & a \le b, \\
p v^{a-b} v^b v^{*b} p = p v^{a-b} p v^b v^{*b} & \mbox{if} & a \ge b,
\end{array} \right.
\]
we obtain by the same argument that $a=b$ if $pwp \neq 0$. Thus, $w = v^a v^{*a}$. The case when $w = v^{*a} v^b v^{*c}$ can be treated analogously. \hfill \qed \\[3mm]
An element $k$ of a $C^*$-algebra $\cA$ is called an {\em element of algebraic rank one} if, for every $a \in \cA$, there is a complex number $\alpha$ such that $kak = \alpha k$.
\begin{prop} \label{p15}
Let $m, \, n \in \sN_v$. Then \\[1mm]
$(a)$ $\pi_n$ is a projection of algebraic rank one in $\cA$. \\
$(b)$ $\pi_m$ and $\pi_n$ are Murray-von Neumann equivalent if and only if $m = n$.
\end{prop}
Analogous assertions hold for $\tilde{\pi}_n$ in place of $\pi_n$. \\[3mm]
{\bf Proof.} $(a)$ Every element of $\cA$ is a limit of linear combinations of words in $v$ and $v^*$. It is thus sufficient to show that, for every word $w$, there is an $\alpha \in \sC$ such that $\pi_n w \pi_n = \alpha \pi_n$. If $\pi_n w \pi_n = 0$, this holds with $\alpha = 0$. If $\pi_n w \pi_n = \pi_n pwp \pi_n \neq 0$, then $w = v^a v^{*a}$ for some $a \ge 0$ by Corollary \ref{c14}. In this case,
\[
\pi_n w \pi_n = \pi_n v^a v^{*a} \pi_n = p (v^n v^{*n} - v^{n+1} (v^*)^{n+1}) v^a v^{*a} (v^n v^{*n} - v^{n+1} (v^*)^{n+1}) p.
\]
From Lemma \ref{l10} we infer that
\[
(v^n v^{*n} - v^{n+1} (v^*)^{n+1}) v^a v^{*a} = \left\{ \begin{array}{lcl}
v^n v^{*n} - v^{n+1} (v^*)^{n+1} & \mbox{if} & a \le n, \\
v^a v^{*a} - v^a v^{*a} = 0      & \mbox{if} & a \ge n+1.
\end{array} \right.
\]
Thus,
\[
\pi_n w \pi_n = \left\{ \begin{array}{lcl}
p (v^n v^{*n} - v^{n+1} (v^*)^{n+1})^2 p = \pi_n & \mbox{if} & a \le n, \\
0                                                & \mbox{if} & a \ge n+1,
\end{array} \right.
\]
i.e. $\alpha = 1$ if $a \le n$ and $\alpha = 0$ in all other cases. \\[1mm]
$(b)$ The projections $\pi_m$ and $\pi_n$ are Murray-von Neumann equivalent if and only if $\pi_m \cA \pi_n \neq \{0\}$. So we have to show that $\pi_m \cA \pi_n = \{0\}$ whenever $m \neq n$. Again it is sufficient to show that $\pi_m w \pi_n = 0$ for every word $w$.

Suppose there is a word $w$ such that $\pi_m w \pi_n = \pi_m p w p \pi_n \neq 0$. Then $w = v^a v^{*a}$ for some $a \ge 0$ by Corollary \ref{c14}. The terms in parentheses in
\[
\pi_m w \pi_n = \pi_m v^a v^{*a} \pi_n = p (v^m v^{*m} - v^{m+1} (v^*)^{m+1}) (v^a v^{*a}) (v^n v^{*n} - v^{n+1} (v^*)^{n+1}) p
\]
commute by Lemma \ref{l10}. Since
\[
(v^m v^{*m} - v^{m+1} (v^*)^{m+1})(v^n v^{*n} - v^{n+1} (v^*)^{n+1}) = 0
\]
for $m \neq n$ we conclude that $\pi_m w \pi_n = 0$, a contradiction. \hfill \qed
\begin{lemma} \label{l16}
$(a)$ If $a > n$ or $b > a$, then $v^b v^{*a} \pi_n = 0$. \\
$(b)$ If $b \le a \le n$, then $v^b v^{*a} \pi_n =  (v^*)^{a-b} \pi_n$.
\end{lemma}
{\bf Proof.} $(a)$ One easily checks that $(v^*)^{n+1}p = 0$, which gives the first assertion. Let $b > a$. Then, since $p$ commutes with $v^k v^{*k}$ and $vp = 0$,
\begin{eqnarray*}
v^b v^{*a} \pi_n & = & v^{b-a} v^a v^{*a} (v^n v^{*n} - v^{n+1} (v^*)^{n+1}) p \\
& = & v^{b-a} p v^a v^{*a} (v^n v^{*n} - v^{n+1} (v^*)^{n+1}) = 0.
\end{eqnarray*}
$(b)$ Applying Lemma \ref{l12} to the terms in inner parentheses in
\[
v^b v^{*a} \pi_n = ((v^b v^{*a} v^n) v^{*n} - (v^b v^{*a} v^{n+1}) (v^*)^{n+1}) p,
\]
one can simplify this expression to
\[
((v^*)^{a-b} v^n v^{*n} - (v^*)^{a-b} v^{n+1}) (v^*)^{n+1}) p = (v^*)^{a-b} \pi_n. \hspace*{20mm}\qed
\]
\begin{coro} \label{c17}
$(a)$ If $w$ is a word in $v, \, v^*$, then $w \pi_n \in \{0, \, \pi_n, \, v^* \pi_n, \, \ldots, v^{*n} \pi_n\}$. \\[1mm]
$(b)$ For every $w \in \cA$, $w \pi_n$ is a linear combination of elements $v^{*i} \pi_n$ with $i \in \{0, \, 1, \, \ldots, n\}$. \\[1mm]
$(c)$ Every element of the ideal $\cC_n$ generated by $\pi_n$ is a linear combination of elements $v^{*i} \pi_n v^j$ with $i, \, j \in \{0, \, 1, \, \ldots n\}$.
\end{coro}
In particular, $\cC_n$ is a finite-dimensional $C^*$-algebra. We are now in a position to describe this algebra exactly.
\begin{prop} \label{p18}
$(a)$  For $n \in \sN_v$, the algebra $\cC_n$ is $^*$-isomorphic to $\sC^{(n+1) \times (n+1)}$. \\[1mm]
$(b)$ $\cC_m \cC_n = \{0\}$ whenever $m \neq n$.
\end{prop}
{\bf Proof.} $(a)$ The elements $e_{ij}^{(n)} := v^{*i} \pi_n v^j$ with $i, \, j \in \{0, \, 1, \, \ldots n\}$ span the algebra $\cC_n$ by Corollary \ref{c17} $(c)$. Thus, the assertion will follow once we have shown that these elements form a system of $(n+1) \times (n+1)$ matrix units in the sense that $(e_{ij}^{(n)})^* = e_{ji}^{(n)}$ and
\begin{equation} \label{e7}
e_{ij}^{(n)} e_{kl}^{(n)} = \delta_{jk} e_{il}^{(n)} \qquad \mbox{for all} \; i, \, j, \, k, \, l \in \{0, \, 1, \, \ldots n\},
\end{equation}
with $\delta_{jk}$ the standard Kronecker delta. The symmetry property is clear. To check (\ref{e7}), first let $j=k$. Then
\[
e_{ij}^{(n)} e_{jl}^{(n)} = v^{*i} \pi_n (v^j v^{*j} \pi_n) v^l = v^{*i} \pi_n^2 v^l = e_{il}^{(n)}
\]
by Lemma \ref{l16} $(b)$. If $j > k$, then
\[
e_{ij}^{(n)} e_{kl}^{(n)} = v^{*i} \pi_n (v^j v^{*k} \pi_n) v^l = 0
\]
by Lemma \ref{l16} $(a)$. Finally, if $j < k$, then
\[
e_{ij}^{(n)} e_{kl}^{(n)} = v^{*i} (\pi_n v^j v^{*k}) \pi_n v^l = v^{*i} (v^k v^{*j} \pi_n)^* \pi_n v^l = 0,
\]
again by Lemma \ref{l16} $(a)$. This proves $(a)$. Assertion $(b)$ follows from Proposition \ref{p15} $(b)$. \hfill \qed \\[3mm]
Given a PPI $v$, we let $\cG_v$ stand for the smallest closed ideal which contains all projections $\pi_n$. If $\sN_v$ is empty, then $\cG_v$ is the zero ideal. Let $\sN_v \neq \emptyset$. The ideal generated by a projection $\pi_n$ with $n \in \sN_v$ is isomorphic to $\sC^{(n+1) \times (n+1)}$ by Proposition \ref{p18}, and if $u, \, w$ are elements of $\cA$  which belong to ideals generated by two different projections $\pi_m$ and $\pi_n$, then $uw = 0$ by Proposition \ref{p15} $(b)$. Hence, $\cG_v$ is then isomorphic to the direct sum of all matrix algebras $\sC^{(n+1) \times (n+1)}$ with $n \in \sN_v$.

If $\cA$ is the {\em universal} $C^*$-algebra generated by a PPI $v$, then $\sN_v$ is the set of {\em all} positive integers. Indeed, the algebra $\cS_{\ge 2}(\eT(C))$ introduced in the introduction is generated by the PPI $v := (P_n V P_n)$, and $\sN_v = \sN$ in this concrete setting.
\begin{coro} \label{c19}
If $\cA$ is the universal $C^*$-algebra generated by a PPI $v$, then $\sN_v = \sN$, and $\cG_v$ is isomorphic to the ideal $\cG_{\ge 2} := \cS_{\ge 2}(\eT(C)) \cap \cG$.
\end{coro}
\section{PPI with $\sN_v = \emptyset$}
Our next goal is to describe the $C^*$-algebra $\cA$ which is generated by a PPI $v$ with $\sN_v = \emptyset$. This condition is evidently satisfied if one of the projections $p = e - v^*v$ and $\tilde{p} = e - vv^*$ is zero, in which cases the algebra generated by the PPI $v$ is well known:
\begin{itemize}
\item
If $p = 0$ and $\tilde{p} = 0$, then $v$ is unitary, and $\cA$ is $^*$-isomorphic to $C(X)$ where $X \subseteq \sT$ is the spectrum of $v$ by the Gelfand-Naimark theorem.
\item
If $p = 0$ and $\tilde{p} \neq 0$, then $v$ is a non-unitary isometry, $\cA$ is $^*$-isomorphic to the Toeplitz algebra $\eT(C)$ by Coburn's theorem, and the isomorphism sends $v$ to the forward shift $V$.
\item
If $p \neq 0$ and $\tilde{p} = 0$, then $v$ is a non-unitary co-isometry, again $\cA$ is $^*$-isomorphic to the Toeplitz algebra $\eT(C)$ by Coburn's theorem, and the isomorphism sends $v$ to the backward shift $V^*$.
\end{itemize}
Thus the only interesting case is when $\sN_v = \emptyset$, but $p\neq 0$ and $\tilde{p} \neq 0$. Let $\cC$ and $\widetilde{\cC}$ denote the smallest closed ideals of $\cA$ which contain the projections $p$ and $\tilde{p}$, respectively. For $i, \, j \ge 0$, set
\[
f_{ij} := v^{*i} p v^j \quad \mbox{and} \quad \tilde{f}_{ij} := v^i \tilde{p} v^{*j}.
\]
\begin{lemma} \label{l20}
If $v$ is a PPI with $\sN_v = \emptyset$, then $(f_{ij})_{i, \, j \ge 0}$ is a (countable) system of matrix units, i.e. $f_{ij}^* = f_{ji}$ and
\begin{equation} \label{e8}
f_{ij} f_{kl} = \delta_{jk} f_{il} \qquad \mbox{for all} \; i, \, j, \, k, \, l \ge 0.
\end{equation}
If one of the $f_{ij}$ is non-zero (e.g. if $f_{00} = p \neq 0$), then all $f_{ij}$ are non-zero.
\end{lemma}
An analogous assertion holds for the family of the $\tilde{f}_{ij}$. \\[3mm]
{\bf Proof.} The symmetry condition is evident, and if $f_{ij} = 0$ then $f_{kl} = f_{ki} f_{ij} f_{jl} = 0$ for all $k, \, l$ by (\ref{e8}). Property (\ref{e8}) on its hand will follow once we have shown that
\begin{equation} \label{e9}
p v^j v^{*k} p = \delta_{jk} p \qquad \mbox{for all} \; j, \, k \ge 0.
\end{equation}
The assertion is evident if $j = k = 0$. If $j > 0$ and $k = 0$, then
\begin{eqnarray*}
p v^j p & = & (e - v^*v) v^j (e - v^*v) \\
& = & v^j - v^*v^{j+1} - v^{j-1} (vv^*v) + v^* v^j (vv^*v) \\
& = & v^j - v^*v^{j+1} - v^{j-1}v + v^* v^j v = 0,
\end{eqnarray*}
and (\ref{e9}) holds. Analogously, (\ref{e9}) holds if $j=0$ and $k > 0$. Finally, let $j, \, k > 0$. The assumption $\sN_v = \emptyset$ ensures that
\begin{equation} \label{e10}
p v^{j-1} \tilde{p} = (e - v^*v) v^{j-1} (e - v v^*) = v^{j-1} - v^* v^j - v^j v^* + v^* v^{j+1} v^* = 0
\end{equation}
for all $j \ge 1$. Employing this identity we find
\begin{eqnarray*}
p v^j v^{*k} p & = & (e - v^*v) \, v^j v^{*k} p \; = \; v^j v^{*k} p - (v^* v^{j+1} v^*) \, (v^*)^{k-1} p \\
& = & v^j v^{*k} p - (v^{j-1} - v^* v^j - v^j v^*) (v^*)^{k-1} p \\
& = &(e - v^*v) v^{j-1} (v^*)^{k-1} p.
\end{eqnarray*}
Thus, $p v^j v^{*k} p = p v^{j-1} (v^*)^{k-1} p$ for $j, \, k \ge 1$. Repeated application of this identity finally leads to one of the cases considered before. \hfill \qed
\begin{prop} \label{p21}
Let $\sN_v = \emptyset$ and $p \neq 0$. \\[1mm]
$(a)$ The ideal $\cC$ of $\cA$ generated by $p$ coincides with the smallest closed subalgebra of $\cA$ which contains all $f_{ij}$ with $i, \, j \ge 0$. \\[1mm]
$(b)$ $\cC$ is $^*$-isomorphic to the ideal of the compact operators on a separable infinite-dimensional Hilbert space. \\[1mm]
Analogous assertions hold for the projection $\tilde{p}$, the algebra $\widetilde{\cC}$, and the $\tilde{f}_{ij}$.
\end{prop}
{\bf Proof.} For a moment, write $\cC^\prime$ for the smallest closed subalgebra of $\cA$ which contains all $f_{ij}$ with $i, \, j \ge 0$. The identities
\[
f_{ij} v = v^{*i} p v^j v = f_{i,j+1}, \qquad vf_{0j} = v p v^j = v (e - v^*v) v^j = 0
\]
and, for $i \ge 1$,
\begin{eqnarray*}
v f_{ij} & = & v v^{*i} (e - v^*v) v^j \; = \; v v^{*i} v^j - (v (v^*)^{i+1} v) v^j \\
& = & v v^{*i} v^j + ((v^*)^{i-1} - (v^*)^i v - v v^{*i}) v^j \\
& = & (v^*)^{i-1} v^j - (v^*)^i v v^j \; = \; (v^*)^{i-1} (e - v^*v) v^j \; = \; f_{i-1,j}
\end{eqnarray*}
(where we used the adjoint of (\ref{e10})) and their adjoints show that $\cC^\prime$ is a closed ideal of $\cA$. Since $p = f_{00}$ we conclude that $\cC \subseteq \cC^\prime$. Conversely, we have $f_{ij} = v^{*i} p v^j \in \cC$ for all $i, \, j \ge 0$ whence the reverse inclusion $\cC^\prime \subseteq \cC$. This settles assertion $(a)$.

For assertion $(b)$, note that every $C^*$-algebra generated by a (countable) system of matrix units (in particular, the algebra $\cC^\prime$) is naturally $^*$-isomorphic to the algebra of the compact operators on a separable infinite-dimensional Hilbert space
(see, e.g., Corollary A.9 in Appendix A2 in \cite{Rae1}). \hfill \qed
\begin{lemma} \label{l22}
If $\sN_v = \emptyset$, then $\cC \cap \widetilde{\cC} = \{0\}$.
\end{lemma}
{\bf Proof.} $\cC$ and $\widetilde{\cC}$ are closed ideals. Thus, $\cC \cap \widetilde{\cC} = \cC \widetilde{\cC}$, and we have to show that $f_{ij} \tilde{f}_{kl} = 0  $ for all $i, \, j, \, k, \, l \ge 0$. Since
\[
f_{ij} \tilde{f}_{kl} = (v^{*i} p v^j) \, (v^k \tilde{p} v^{*l}) = v^{*i} \, (p v^{j+k} \tilde{p}) \, v^{*l},
\]
this is a consequence of $\sN_v = \emptyset$. \hfill \qed \\[3mm]
Remember that $p \neq 0$ and $\tilde{p} \neq 0$. From the preceding lemma we conclude that the mapping
\[
\cA \to \cA/\cC \times \cA/\widetilde{\cC}, \quad w \mapsto (w + \cC, \, w + \widetilde{\cC})
\]
is an injective $^*$-homomorphism; thus $\cA$ is $^*$-isomorphic to the $C^*$-subalgebra of $\cA/\cC \times \cA/\widetilde{\cC}$ generated by $(v + \cC, \, v + \widetilde{\cC})$. The element $v + \cC$ is an isometry in $\cA/\cC$ (since $e - v^*v \in \cC$), but it is not unitary (otherwise $e - vv^* \in \tilde{\cC}$ would be a non-zero element of $\cC$, in contradiction with Lemma \ref{l22}). Analogously, $v + \widetilde{\cC}$ is a non-unitary co-isometry. By Coburn's Theorem, there are $^*$-isomorphisms $\mu : \cA/\cC \to \eT(C)$ and $\tilde{\mu} : \cA/\tilde{\cC} \to \eT(C)$ which map $v + \cC \mapsto V$ and $v + \tilde{\cC} \mapsto V^*$, respectively. But then
\[
\mu \times \tilde{\mu} : \cA/\cC \times \cA/\widetilde{\cC} \to \eT(C) \times \eT(C), \quad (a, \, \tilde{a}) \mapsto (\mu(a), \, \tilde{\mu} (\tilde{a}))
\]
is a $^*$-isomorphism which maps the $C^*$-subalgebra of $\cA/\cC \times \cA/\widetilde{\cC}$ generated by $(v + \cC, \, v + \widetilde{\cC})$ to the $C^*$-subalgebra of $\eT(C) \times \eT(C)$ generated by the pair $(V, \, V^*$). The latter algebra has been identified in Corollary \ref{c4}. Summarizing we get:
\begin{prop} \label{p23}
Let the $C^*$-algebra $\cA$ be generated by a PPI $v$ with $\sN_v = \emptyset$ and $p \neq 0$ and $\tilde{p} \neq 0$. Then $\cA$ is $^*$-isomorphic to the algebra $\cS(\eT(C))/\cG$
$($likewise, to $\cS_{\ge 2}(\eT(C))/\cG_{\ge 2})$, and the isomorphism sends $v$ to $(P_n V P_N)_{n \ge 1} + \cG$ $($likewise, to $(P_n V P_N)_{n \ge 2} + \cG_{\ge 2})$.
\end{prop}
\section{The general case}
We are now going to finish the proof of Theorem \ref{t1}. For we think of $\cA$ as being faithfully represented as a $C^*$-algebra of bounded linear operators on a separable infinite-dimensional Hilbert space $H$ (note that $\cA$ is finitely generated, hence separable). As follows easily from (\ref{e7}), $z_n := \sum_{i=0}^n e_{ii}^{(n)}$ is the identity element of $\cC_n$. So we can think of the $z_n$ as orthogonal projections on $H$. Moreover, these projections are pairwise orthogonal by Proposition \ref{p18} $(b)$. Thus, the operators $P_n := \sum_{i=1}^n z_n$ form an increasing sequence of orthogonal projections on $H$. Let $P \in L(H)$ denote the least upper bound of that sequence (which then is the limit of the $P_n$ in the strong operator topology). Clearly, $P$ is an orthogonal projection again (but note that $P$ does not belong to $\cA$ in general).
\begin{lemma} \label{l24}
$(a)$ Every $z_n$ is a central projection of $\cA$. \\[1mm]
$(b)$ $P$ commutes with every element of $\cA$.
\end{lemma}
{\bf Proof.} Assertion $(b)$ is a consequence of $(a)$. We show that
\begin{eqnarray*}
z_n & = & \sum_{i=0}^n v^{*i} \pi_n v^i = \sum_{i=0}^n v^{*i} p (v^n v^{*n} - v^{n+1} (v^*)^{n+1}) v^i \\
& = & \sum_{i=0}^n v^{*i} (e - v^*v) (v^n v^{*n} - v^{n+1} (v^*)^{n+1}) v^i
\end{eqnarray*}
commutes with $v$. Indeed,
\begin{eqnarray*}
v z_n & = & v (e - v^*v) (v^n v^{*n} - v^{n+1} (v^*)^{n+1}) \\
&& \qquad + \sum_{i=1}^n v v^{*i} (e - v^*v) (v^n v^{*n} - v^{n+1} (v^*)^{n+1}) v^i \\
& = & \sum_{i=1}^n v v^{*i} (e - v^*v) v^n (v^{*n} - v (v^*)^{n+1}) v^i \\
& = & \sum_{i=1}^n (v v^{*i} v^n - v (v^*)^{i+1} v^{n+1}) (v^{*n} - v (v^*)^{n+1}) v^i \\
& = & \sum_{i=1}^n ((v^*)^{i-1} v^n - (v^*)^{i} v^{n+1}) (v^{*n} - v (v^*)^{n+1}) v^i \quad \mbox{(by Lemma \ref{l12})} \\
& = & \sum_{i=1}^n (v^*)^{i-1} (e - v^* v) (v^n v^{*n} - v^{n+1} (v^*)^{n+1}) v^i \\
& = & \sum_{i=0}^{n-1} v^{*i} (e - v^* v) (v^n v^{*n} - v^{n+1} (v^*)^{n+1}) v^i v \\
& = & \sum_{i=0}^n v^{*i} (e - v^* v) (v^n v^{*n} - v^{n+1} (v^*)^{n+1}) v^i v \\
&& \qquad - (v^*)^{n} (e - v^* v) (v^n v^{*n} - v^{n+1} (v^*)^{n+1}) v^{n+1} \\
& = & \sum_{i=0}^n v^{*i} (e - v^* v) (v^n v^{*n} - v^{n+1} (v^*)^{n+1}) v^i v \; = \; z_n v
\end{eqnarray*}
again by Lemma \ref{l12}. Thus, $v z_n = z_n v$. Since $z_n = z_n^*$, this implies that $z_n$ also commutes with $v^*$ and, hence, with every element of $\cA$. \hfill \qed \\[3mm]
Consequently, $\cA = P \cA P \oplus (I-P) \cA (I-P)$ where $I$ stands for the identity operator on $H$. We consider the summands of this decomposition separately. The part $(I-P) \cA (I-P)$ is generated by the PPI $v^\prime := (I-P)v(I-P)$. Since
\begin{eqnarray*}
(I-P) p v^n \tilde{p} v^{*n} p (I-P) & = & (I-P) \pi_n (I-P) \\
& = & (I-P) z_n e_{00}^{(n)} \pi_n (I-P) = 0,
\end{eqnarray*}
we conclude that $\sN_{v^\prime} = \emptyset$. Thus, this part of $\cA$ is described by Proposition \ref{p23}.

The part $P \cA P$ is generated by the PPI $PvP$. It follows from the definition of $P$ that $\sN_{PvP} = \sN_v$ and that $\cG_{PvP} = P\cG_vP = \cG_v$. We let $\prod_{n \in \sN_v} \cC_n$ stand for the direct product of the algebras $\cC_n$ and consider the mapping
\begin{equation} \label{e11}
P \cA P \to \prod_{n \in \sN_v} \cC_n, \quad PAP \mapsto (z_n PAP z_n)_{n \in \sN_v} = (z_n A z_n) _{n \in \sN_v}.
\end{equation}
If $z_n A z_n = 0$ for every $n \in \sN_v$, then
\[
PAP = \sum_{m, n \in \sN_v} z_m PAP z_n = \sum_{n \in \sN_v} z_n A z_n = 0.
\]
Thus, the mapping (\ref{e11}) is injective, and the algebra $P \cA P$ is $^*$-isomorphic to the $C^*$-subalgebra of $\prod_{n \in \sN_v} \cC_n$ generated by the sequence $(z_n v z_n) _{n \in \sN_v}$. Further we infer from Proposition \ref{p18} $(a)$ that $\cC_n$ is isomorphic to $\sC^{(n+1) \times (n+1)}$ if $n \in \sN_v$. We are going to make the latter isomorphism explicit. For we note that
\begin{eqnarray*}
e_{ii}^{(n)} v e_{jj}^{(n)} & = & v^{*i} \pi_n v^{i+1} v^{*j} \pi_n v^j \\
& = &  \left\{ \begin{array}{lcl}
v^{*i} \pi_n v^{i+1} (v^*)^{i+1} \pi_n v^{i+1} & \mbox{if} & i+1 = j, \\
0                                              & \mbox{if} & i+1 \neq j
\end{array} \right.
\quad \mbox{(by Corollary \ref{c14})} \\
& = &  \left\{ \begin{array}{lcl}
v^{*i} \pi_n v^{i+1} & \mbox{if} & i+1 = j, \\
0                    & \mbox{if} & i+1 \neq j
\end{array} \right.
\quad \mbox{(by Lemma \ref{l16})} \\
& = &  \left\{ \begin{array}{lcl}
e_{i,i+1}^{(n)}& \mbox{if} & i+1 = j, \\
0           & \mbox{if} & i+1 \neq j.
\end{array} \right.
\end{eqnarray*} 
We choose a unit vector $e_i^{(n)}$ in the range of $e_{ii}^{(n)}$ (recall Proposition \ref{p15} $(a)$), and let $f_i^{(n)}$ stand for the $n+1$-tuple $(0, \ldots, 0, \, 1, \, 0, \, \ldots, \, 0)$ with the 1 at the $i$th position. Then $(e_i^{(n)})_{i=0}^n$ forms an orthonormal basis of $\im z_n$, $(f_i^{(n)})_{i=0}^n$ forms an orthonormal basis of $\sC^{n+1}$, the mapping $e_i^{(n)} \mapsto f_{n-i}^{(n)}$ extends to a linear bijection from $\im z_n$ onto $\sC^{n+1}$, which finally induces a $^*$-isomorphism $\xi_n$ from $\cC_n \cong L(\im z_n)$ onto $\sC^{(n+1) \times (n+1)} \cong L(\sC^{n+1})$. Then 
\[
\xi : \prod_{n \in \sN_v} \cC_n \to \prod_{n \in \sN_v} \sC^{(n+1) \times (n+1)}, \quad (A_n) \mapsto (\xi_n(A_n))
\]
is a $^*$-isomorphism which maps the $C^*$-subalgebra of $\prod_{n \in \sN_v} \cC_n$ generated by the sequence $(z_n v z_n) _{n \in \sN_v}$ to the $C^*$-subalgebra of $\prod_{n \in \sN_v} \sC^{(n+1) \times (n+1)}$ generated by the sequence $(V_{n+1}) _{n \in \sN_v}$, where $V_n$ is the matrix described in Example $(c)$. Note that $V_n$ is just the $n \times n$th finite section $P_n V P_n$ of the forward shift operator.

If now $\cA$ is the universal algebra generated by a PPI $v$, then $\sN_v = \sN$, as we observed in Corollary \ref{c19}. Thus, in this case, the algebra $P \cA P$ is $^*$-isomorphic to the smallest $C^*$-subalgebra of $\cF = \prod_{n \ge 1} \sC^{n \times n}$ generated by the sequence $(P_n V P_n)$, i.e. to the $C^*$-algebra $\cS_{\ge 2}(\eT(C))$.

It remains to explain what happens with the part $(I-P) \cA (I-P)$ of $\cA$. The point is that the quotient $P \cA P/P \cG_v P$ is generated by a PPI $u$ for which  
$\sN_u$ is empty. We have seen in Proposition \ref{p23} that both this quotient and the algebra $(I-P) \cA (I-P)$ are canonically $^*$-isomorphic to $\cS_{\ge 2}(\eT(C))/\cG_{\ge 2}$. Thus, there is a $^*$-homomorphism from $P \cA P$ onto $(I-P) \cA (I-P)$ which maps the generating PPI $PvP$ of $P \cA P$ to the generating PPI $(I-P) v (I-P)$ of $(I-P) \cA (I-P)$. Hence, if $\cA$ is the universal $C^*$-algebra generated by a PPI, then already $P \cA P$ has the universal property, and $\cA \cong P \cA P \cong \cS_{\ge 2}(\eT(C))$. \hfill \qed
{\small Author's address: \\[3mm]
Steffen Roch, Technische Universit\"at Darmstadt, Fachbereich Mathematik, Schlossgartenstrasse 7, 64289 Darmstadt, Germany. \\
E-mail: roch@mathematik.tu-darmstadt.de}
\end{document}